\documentclass[11pt]{amsart}

\usepackage{amsmath}
\usepackage{amssymb}
\usepackage{mathrsfs}

\newtheorem{theorem}{Theorem}[section]
\newtheorem{claim}[theorem]{Claim}

\newtheorem{lemma}[theorem]{Lemma}

\newtheorem{corollary}[theorem]{Corollary}

\theoremstyle{definition}
\newtheorem{definition}[theorem]{Definition}

\theoremstyle{remark}
\newtheorem{remark}[theorem]{Remark}

\newcount\skewfactor
\def\mathunderaccent#1#2 {\let\theaccent#1\skewfactor#2
\mathpalette\putaccentunder}
\def\putaccentunder#1#2{\oalign{$#1#2$\crcr\hidewidth
\vbox to.2ex{\hbox{$#1\skew\skewfactor\theaccent{}$}\vss}\hidewidth}}

\def\smallbox#1{\leavevmode\thinspace\hbox{\vrule\vtop{\vbox
   {\hrule\kern1pt\hbox{\vphantom{\tt/}\thinspace{\tt#1}\thinspace}}
   \kern1pt\hrule}\vrule}\thinspace}

\newcommand{\bool}{{\bf B}}

\newcommand{\cf}{{\rm cf}}
\newcommand{\Depth}{{\rm Depth}}
\newcommand{\Length}{{\rm Length}}

\def\qedref#1{$\qed_{\reforiginal{#1}}$}

\title{Depth$^+$ and Length$^+$ of Boolean Algebras}
\author{Shimon Garti}
\address{Einstein Institute of Mathematics,
 The Hebrew University of Jerusalem,
 Jerusalem 91904, and Ben-Gurion University of the Negev, Beer-Sheva, Israel}
\email{shimon.garty@mail.huji.ac.il}
\thanks{}

\author{Saharon Shelah}
\address{Institute of Mathematics
 The Hebrew University of Jerusalem
 Jerusalem 91904, Israel
 and  Department of Mathematics
 Rutgers University
 New Brunswick, NJ 08854, USA}
\email{shelah@math.huji.ac.il}
\urladdr{http://www.math.rutgers.edu/\char`\~shelah}
\thanks{Research supported by ERC grant 338821. Publication 974 of the second author}


\begin{document}
\let\labeloriginal\label
\let\reforiginal\ref
\def\ref#1{\reforiginal{#1}}
\def\label#1{\labeloriginal{#1}}

\begin{abstract}
Suppose that $\kappa=\cf(\kappa), \lambda>\cf(\lambda)=\kappa^+$ and $\lambda=\lambda^\kappa$. We prove that there exist a sequence $\langle\bool_i:i<\kappa\rangle$ of Boolean algebras and an ultrafilter $D$ over $\kappa$ so that $\lambda=\prod\limits_{i<\kappa}\Depth^+(\bool_i)/D< \Depth^+(\prod\limits_{i<\kappa}\bool_i/D)=\lambda^+$. An identical result holds also for $\Length^+$.
The proof is carried in ZFC, and it holds even above large cardinals.
\end{abstract}

\maketitle

\newpage

\section{Introduction}

The monograph of Monk, \cite{MR1393943}, lists many cardinal invariants on Boolean algebras. One of them is called $\Depth$, and it concerns with well ordered subsets of Boolean algebras. But there are two variations of this invariant, as can be seen from the following:

\begin{definition}
\label{dddplus}
$\Depth$ and $\Depth^+$ of Boolean algebras. \newline 
Let $\bool$ be a Boolean Algebra.
\begin{enumerate}
\item [$(\aleph)$] $\Depth(\bool)= \sup\{\theta : \exists \bar b=(b_\gamma:\gamma<\theta), \hbox{ increasing\ sequence\ in } \bool\}$.
\item [$(\beth)$] $\Depth^+ (\bool)=\sup \{\theta^+:\exists \bar b=(b_\gamma:\gamma<\theta), \hbox{ increasing\ sequence\ in } \bool\}$.
\end{enumerate}
\end{definition}

Another invariant is the $\Length$. Again, we have two variations:

\begin{definition}
\label{lllplus}
$\Length$ and $\Length^+$ of Boolean algebras. \newline 
Let $\bool$ be a Boolean Algebra.
\begin{enumerate}
\item [$(\aleph)$] $\Length(\bool)= \sup\{\theta : \exists A\subseteq\bool, |A|=\theta$ such that $A$ is linearly-ordered by $<_\bool\}$.
\item [$(\beth)$] $\Length^+(\bool)= \sup\{\theta^+ : \exists A\subseteq\bool, |A|=\theta$ such that $A$ is linearly-ordered by $<_\bool\}$.
\end{enumerate}
\end{definition}

Take a look at the definitions of $\Depth$ and $\Depth^+$.
At first glance it seems that the difference between these two variants has a technical nature. The theme of this paper is to show that the difference is important, and the `correct' definition should be $\Depth^+$. 

Let us consider a Boolean algebra $\bool$, such that $\Depth(\bool)$ is a limit cardinal $\lambda$. It might happen that $\lambda$ is not attained (i.e., there is a chain of length $\theta$ for every $\theta<\lambda$ in $\bool$, but no chain of length $\lambda$), and it might happen that $\lambda$ is attained (i.e., there is a chain of length $\lambda$ in $\bool$). In both cases, $\Depth(\bool)=\lambda$. On the other hand, $\Depth^+(\bool)=\lambda$ in the first scene, but $\Depth^+(\bool)=\lambda^+$ in the second. The conclusion is that $\Depth$ is less informative than $\Depth^+$.

The little example above is very simple, but the same phenomenon reflects in other related problems, including the problem of ultraproducts. In this paper we deal with this construction. Let us try to sketch the background and history of the problem.

Suppose $inv$ is any cardinal invariant on Boolean algebras. Given a sequence $\langle\bool_i:i<\kappa\rangle$ of Boolean algebras and an ultrafilter $D$ on $\kappa$, we can walk in two courses. In the algebraic route we define a new Boolean algebra $\bool=\prod\limits_{i<\kappa}\bool_i/D$. Having the algebra $\bool$, we compute $inv(\bool)$. In the set theoretical route we produce a sequence of cardinals, $\langle inv(\bool_i):i<\kappa\rangle$, say $\theta_i=inv(\bool_i)$ for every $i<\kappa$. Now we compute $\prod\limits_{i<\kappa}\theta_i/D$.

Monk investigates systematically the relationship between these two routes. We are looking for constructions which give strict inequalities (in both directions). We are also interested in the consistency power of these constructions. The most basic problem here is if such a construction can be carried out in ZFC.

It is consistent with ZFC that $\prod\limits_{i<\kappa}\Depth(\bool_i)/D\leq \Depth(\prod\limits_{i<\kappa}\bool_i/D)$ for every ultrafilter $D$ and every sequence $\langle\bool_i:i<\kappa\rangle$ (see theorem 4.14 in \cite{MR1393943}), hence no ZFC counterexample is available. But what about a ZFC example of the strict relation $\prod\limits_{i<\kappa}\Depth(\bool_i)/D <\Depth(\prod\limits_{i<\kappa}\bool_i/D)$? This question is problem number 12 in Monk's list. A parallel problem arises for the $\Length$ invariant (this is labeled as problem number 22 in the same list).

There is a meaningful difference between these problems. Problem number 12 is still open, and we have some restrictions on the (tentative) existence of a ZFC construction which gives $\prod\limits_{i<\kappa}\Depth(\bool_i)/D <\Depth(\prod\limits_{i<\kappa}\bool_i/D)$. First, if $\lambda>\cf(\lambda)=\aleph_0$ then such an example is ruled out (see \cite{MR2822491} and \cite{MR2916079}). Second, the discrepancy (if exists at all) is limited to one cardinal (under the assumption $\lambda^\kappa=\lambda$, see \cite{MR2386531}).

Problem number 22 (about $\Length$) has been solved (in \cite{MR1728851}, Theorem 15.14). The gap in \cite{MR1728851} is one cardinal, but it seems that a larger gap is possible (and we hope to prove it in a subsequent work). Likewise, strict inequalities for $\Length$ were forced in \cite{MR1674385} under some large cardinals assumptions before the ZFC theorem has been discovered. So our knowledge about $\Length$ is deeper than our knowledge about $\Depth$ (with respect to ultraproducts).

Anyway, using the more informative definitions of $\Depth^+$ and $\Length^+$ yields a plenty of ZFC counterexamples, as we shall try to prove in the present work. We also direct the reader to \cite[\S 4]{MR3580893} in which related results are proved.

Our notation is standard. We follow the terminology of \cite{MR991565} and \cite{MR1393943} in general. We shall use the notion of a regular ultrafilter, so we need the following definition:

\begin{definition}
\label{rrrr}
Regular Ultrafilters. \newline 
Let $D$ be an ultrafilter on $\kappa$. \newline 
$D$ is regular if there exists a sequence $\langle W_i:i<\kappa\rangle$, each $W_i$ belongs to $[\kappa]^{<\aleph_0}$, and $\{i<\kappa:\zeta\in W_i\}\in D$ for every $\zeta\in\kappa$.
\end{definition}

The property of regular ultrafilters to be used in the main theorem is that $\prod\limits_{i<\kappa}\lambda_i/D=\lambda^\kappa$, in particular it equals $\lambda$ if we choose a cardinal which satisfies $\lambda^\kappa=\lambda$ as in the theorem below.
This concept and basic fact goes back to Keisler, see \cite{MR0409165}.

We shall make use of the Delta-system lemma. For the general theorem and proof, one may consult \cite{MR597342}.
We need just the simplest form which says that if $\theta_\varepsilon$ is an uncountable regular cardinal and $F_\varepsilon$ is a collection of $\theta_\varepsilon$-many finite sets, then there exists a finite set $r_\varepsilon$ and $I_\varepsilon\in[F_\varepsilon]^{\theta_\varepsilon}$ so that $\{x,y\}\in [I_\varepsilon]^2\Rightarrow x\cap y=r_\varepsilon$. By abuse of notation, we may assume that $r_\varepsilon$ is a set of natural numbers which are the indices of the members in the finite sets of $I_\varepsilon$. For a club set $E$ let ${\rm acc}(E)$ be the set of accumulation points of $E$, i.e. the set $\{\delta:\delta={\rm sup}(E\cap\delta)\}$.

Dealing with Boolean algebras, we quote a specific case of Sikorski's extension theorem. A detailed proof can be found in \cite{MR991565}:

\begin{theorem}
\label{ssikorski} Extending homomorphisms. \newline 
Let $\bool_1$ be a Boolean algebra, generated freely by $\langle x_\gamma:\gamma<\mu\rangle$ except some set $\Gamma\subseteq \{(x_\alpha\leq x_\beta) :\alpha,\beta<\mu\}$ of relations between the generators.
Assume $\bool_2$ is another Boolean algebra, and a function $f$ is defined on $\langle x_\gamma:\gamma<\mu\rangle$ into $\bool_2$ such that $(x_\alpha\leq x_\beta) \in\Gamma\Rightarrow f(x_\alpha)\leq_{\bool_2}f(x_\beta)$. \newline 
Then there is a homomorphism $\hat{f}:\bool_1\rightarrow\bool_2$ which extends $f$.
\end{theorem} 

\hfill \qedref{ssikorski}

Assume $0<\gamma_i<\kappa^+$ for every $i<\kappa$. One can form the product $\prod\limits_{i<\kappa}\gamma_i/D$, when $D$ is an ultrafilter on $\kappa$. Each member of the product is an equivalence class of functions in $\prod\limits_{i<\kappa}\gamma_i$. The equivalence relation is defined by $D$, i.e., $f<_Dg\Leftrightarrow\{i<\kappa: f(i)<g(i)\}\in D$. The following is known (see \cite{MR1393943} p. 90, and \cite{MR1083551} chapter VI \S 3):

\begin{theorem}
\label{xxxx} Increasing chains in ultraproducts. 
\begin{enumerate}
\item [$(a)$] For every $\kappa\geq\aleph_0$ and every uniform ultrafilter $D$ on $\kappa$ there exists an increasing chain of length $\kappa^+$ in $\kappa^\kappa/D$.
\item [$(b)$] For every $\kappa\geq\aleph_0$ and $\kappa^+\leq\partial= \cf(\partial)\leq 2^\kappa$ there exists a regular ultrafilter $D$ on $\kappa$ so that in $\kappa^\kappa/D$ (and even in $\omega^\kappa/D$) there exists an increasing chain of length $\partial$.
\end{enumerate}
\end{theorem}

\hfill \qedref{xxxx}

\newpage

\section{Ultraproducts of Boolean algebras}

Let us begin with the following lemma:

\begin{lemma}
\label{tttransitivity} The transitivity lemma. \newline 
Suppose $D$ is a uniform ultrafilter on $\kappa, \partial\in[\kappa^+, 2^\kappa], 0<\gamma_i<\kappa^+$ for every $i<\kappa$ and in $\prod\limits_{i<\kappa}\gamma_i/D$ there exists an increasing chain of length $\partial$. \newline 
Then one can choose $\langle <_i:i<\kappa\rangle$ and $\langle g_i:i<\kappa\rangle$ so that:
\begin{enumerate}
\item [$(a)$] $<_i$ is a partial order on $\partial$, for every $i<\kappa$.
\item [$(b)$] $g_i$ is a function from $\partial$ into $\gamma_i$, satisfies $\zeta<_i\varepsilon \Rightarrow g_i(\zeta)<g_i(\varepsilon)$ for every $i<\kappa$.
\item [$(c)$] If $\zeta<\varepsilon<\partial$ then $\{i<\kappa: \zeta<_i \varepsilon\}\in D$.
\end{enumerate}
\end{lemma}

\par \noindent \emph{Proof}. \newline 
Fix any sequence $\langle f_\varepsilon:\varepsilon<\partial\rangle$ so that $f_\varepsilon\in \prod\limits_{i<\kappa}\gamma_i/D$ for every $\varepsilon<\partial$, and $\zeta<\varepsilon<\partial \Rightarrow f_\zeta <_D f_\varepsilon$. Such a sequence exists by the assumptions of the lemma. For every $i<\kappa$, define:

$$
\zeta<_i\varepsilon\Leftrightarrow (\zeta<\varepsilon<\partial) \wedge (f_\zeta(i) < f_\varepsilon(i)).
$$

As $<_i$ is a partial order over $\partial$ for every $i<\kappa$, part $(a)$ is satisfied. For every $\varepsilon<\partial$ we define $g_i(\varepsilon)=f_\varepsilon(i)$, so $g_i:\partial\rightarrow\gamma_i$ for every $i<\kappa$. Notice that $\zeta<_i\varepsilon$ implies $g_i(\zeta)=f_\zeta(i)< f_\varepsilon(i)=g_i(\varepsilon)$, hence part $(b)$ is satisfied as well. Finally, if $\zeta<\varepsilon<\partial$ then $f_\zeta <_D f_\varepsilon$ which ammounts to $\{i<\kappa:\zeta<_i\varepsilon\}= \{i<\kappa:f_\zeta(i)<f_\varepsilon(i)\}\in D$, hence part $(c)$ is established and the proof is accomplished.

\hfill \qedref{tttransitivity}

This lemma enables us to define our Boolean algebras in the main theorem. We shall use the lemma in order to make sure that the order of the Boolean algebras is transitive. We need another lemma, which says that a special kind of a Delta-system can be created on a singular cardinal $\lambda$ with uncountable cofinality:

\begin{lemma}
\label{dddd} The singular Delta-system. \newline 
Suppose $\lambda>\cf(\lambda)=\partial>\aleph_0$, and $\{u_\alpha:\alpha<\lambda\}$ is a collection of finite sets. Assume $\langle\theta_\varepsilon:\varepsilon<\partial\rangle$ is an increasing continuous sequence of cardinals which tends to $\lambda$ so that $\theta_0=0, \theta_1>\partial$ and $\theta_{\varepsilon+1}$ is a regular cardinal for every $\varepsilon<\partial$. \newline 
There is a set $B\in[\lambda]^\lambda$ and an unbounded subset $T\in[\partial]^\partial$ such that for every $\gamma_0,\gamma_1\in B, \gamma_0<\gamma_1$ we have the following:
\begin{enumerate}
\item [$(a)$] If $\gamma_0,\gamma_1\in[\theta_\varepsilon, \theta_{\varepsilon+1})$ for some $\varepsilon\in T$, then $u_{\gamma_0}\cap u_{\gamma_1}=r_\varepsilon$ for some fixed finite set $r_\varepsilon$.
\item [$(b)$] If $\gamma_0\in[\theta_\varepsilon, \theta_{\varepsilon+1})$, $\gamma_1\in[\theta_\zeta, \theta_{\zeta+1})$ and $\varepsilon<\zeta$ are from $T$, then $u_{\gamma_0}\cap u_{\gamma_1}=r_*$ for some fixed finite set $r_*$.
\item [$(c)$] $r_{\varepsilon_0}\cap r_{\varepsilon_1}=r_*$ for every $\varepsilon_0<\varepsilon_1$ from $T$.
\item [$(d)$] $|B\cap[\theta_\varepsilon, \theta_{\varepsilon+1})|=\theta_{\varepsilon+1}$, for every $\varepsilon\in T$.
\end{enumerate}
\end{lemma}

\par\noindent\emph{Proof}.\newline 
For every $\varepsilon<\partial$ we have $\theta_{\varepsilon+1}$-many members in the collection $\{u_\alpha:\alpha\in [\theta_\varepsilon, \theta_{\varepsilon+1})\}$, hence there exists $I_\varepsilon\subseteq [\theta_\varepsilon, \theta_{\varepsilon+1})$, $|I_\varepsilon|=\theta_{\varepsilon+1}$ and a fixed finite root $r_\varepsilon$ so that:

$$
\gamma_0,\gamma_1\in I_\varepsilon, \gamma_0<\gamma_1 \Rightarrow u_{\gamma_0}\cap u_{\gamma_1}=r_\varepsilon.
$$

Notice that $I_\varepsilon$ satisfies part $(a)$ of the lemma, and consequently every shrinking of $I_\varepsilon$ satisfies it. Since $\partial>\aleph_0$ there exists a set $T\in[\partial]^\partial$ such that $\{r_\varepsilon:\varepsilon\in T\}$ is a Delta-system, and $r_*$ is the root. This gives us part $(c)$ of the lemma.

For every $\varepsilon\in T$ let $I^-_\varepsilon$ be the following set:

$$
\{\gamma\in I_\varepsilon: [(u_\gamma\setminus r_\varepsilon )\cap\bigcup\limits_{\varepsilon\in T} r_\varepsilon\neq \emptyset] \bigvee [(u_\gamma\setminus r_\varepsilon )\cap\bigcup\limits_{\beta<\theta_\varepsilon} u_\beta\neq \emptyset]\}.
$$

Clearly, $|I_\varepsilon^-|\leq\theta_\varepsilon+\partial$ for every $\varepsilon\in T$. Consequently, $|I_\varepsilon\setminus  I_\varepsilon^-|=\theta_{\varepsilon+1}$ for every $\varepsilon\in T$, hence $B=\bigcup \{I_\varepsilon\setminus  I_\varepsilon^-: \varepsilon\in T\}$ is a member of $[\lambda]^\lambda$. We claim that $B$ is as required.

Indeed, part $(a)$ holds for every $I_\varepsilon$, so also for $I_\varepsilon\setminus  I_\varepsilon^-$. Part $(c)$ has been established, and part $(d)$ follows from the equality $|I_\varepsilon\setminus  I_\varepsilon^-|= \theta_{\varepsilon+1}$. Part $(b)$ follows from removing $I_\varepsilon^-$ (at each $\varepsilon\in T$) which gives $r_*$ as the intersection of every pair of members from distinct layers.

\hfill \qedref{dddd}

\begin{remark}
\label{tttt} A parallel statement can be phrased upon replacing the finite sets $u_\alpha$ by finite sequences $\bar{\gamma}_\alpha$. We shall use, below, the sequence version (the proof is the same, but the notation is more cumbersome).
\end{remark}

\hfill \qedref{tttt}

We can state now the main result of the paper:

\begin{theorem}
\label{mt} The main theorem. \newline 
Assume $\lambda>\cf(\lambda)=\partial, \partial\in[\kappa^+,2^\kappa]$ and $D$ is a uniform ultrafilter on $\kappa$ which satisfies the conclusion of Theorem \ref{xxxx}. \newline 
Then one can find $\langle\bool_i:i<\kappa\rangle$ such that:
\begin{enumerate}
\item [$(\aleph)$] $\Depth^+(\bool_i)\leq\lambda$ for every $i<\kappa$.
\item [$(\beth)$] $\Depth^+(\bool)\geq\lambda^+$ and equality holds if $\lambda^\kappa=\lambda$.
\end{enumerate}
Consequently, $\prod\limits_{i<\kappa}\Depth^+(\bool_i)/D< \Depth^+(\prod\limits_{i<\kappa}\bool_i/D)$.
\end{theorem}

The idea is to define Boolean algebras which are `free enough' to supply many homomorphisms on each Boolean algebra. We shall create this algebra such that if $\langle b_\gamma:\gamma<\lambda\rangle$ is an increasing chain then one can find two members $b_{\gamma_1}<b_{\gamma_2}$ and designate $f:\bool_i\rightarrow\bool_i$ so that $f(b_{\gamma_1})=b_{\gamma_2}$ and $f(b_{\gamma_2})=b_{\gamma_1}$. The existence of this homomorphism is based on the fact that the length of the chain is $\lambda$. This yields a contradiction, since homomorphism (in Boolean algebras) is order preserving. Consequently, we know that no increasing chains of length $\lambda$ exist in $\bool_i$ for every $i<\kappa$, hence part $(\aleph)$ holds. On the other hand, using Lemma \ref{tttransitivity} (c) for our ultrafilter, we will be able to introduce a $\lambda$-chain in the product algebra.

\medskip

\par\noindent\emph{Proof}. \newline 
Assume there are $\kappa,\partial,\lambda$ as in the assumptions of the theorem (notice that for every infinite cardinal $\kappa$, the cardinal $\lambda=\beth_{\partial}(\aleph_0)$ can serve; similarly $\beth_{\delta}(\aleph_0)$ for any ordinal $\delta$ of cofinality $\partial$). Let $D$ be a uniform ultrafilter on $\kappa$ which satisfies the demands in Theorem \ref{xxxx}(b). Note that $D$ can be chosen as a regular ultrafilter. 
Let $\langle \theta_\varepsilon: \varepsilon<\partial\rangle$ be an increasing continuous sequence of cardinals which tends to $\lambda$ such that $\theta_0=0, \theta_1>\partial$, and each $\theta_{\varepsilon+1}$ is regular.

Let $\xi(\alpha)$ be ${\rm min}\{\varepsilon:\theta_\varepsilon\leq \alpha<\theta_{\varepsilon+1}\}$ for every $\alpha<\lambda$. $\xi$ is a `block' function, and $\xi(\alpha)$ determines the unique interval $[\theta_\varepsilon,\theta_{\varepsilon+1})$ which $\alpha$ belongs to.
For every $i<\kappa$ set $\Gamma_i=\{(x^i_\alpha<x^i_\beta):[\alpha<\beta \wedge \xi(\alpha)=\xi(\beta)] \bigvee [\xi(\alpha) <_i \xi(\beta)]\}$.
We define $\bool_i$ as the Boolean algebra generated freely from $\{x_\alpha^i:\alpha<\lambda\}$, except the relations in $\Gamma_i$. Lemma \ref{tttransitivity} tells us that $\bool_i$ is a Boolean algebra.
This definition accomplishes the construction of the Boolean algebras, and recall that $\bool$ is the ultraproduct algebra.

We shall elicit an increasing sequence $\langle y_\gamma:\gamma<\lambda\rangle$ of members of $\bool$. For every $\gamma<\lambda$ we set $y_\gamma=\langle x^j_\gamma:j<\kappa\rangle/D$.
Suppose $\gamma_0<\gamma_1<\lambda$, so $\xi(\gamma_0)\leq\xi(\gamma_1)$. If $\xi(\gamma_0)=\xi(\gamma_1)$ then for every $i<\kappa$ we have $x^i_{\gamma_0}<_{\bool_i} x^i_{\gamma_1}$, and since $\kappa\in D$ we conclude that $y_{\gamma_0}<_\bool y_{\gamma_1}$. 
If $\xi(\gamma_0)< \xi(\gamma_1)$ then $\{i<\kappa:\xi(\gamma_0)<_i \xi(\gamma_1)\}\in D$ and consequently $\{i<\kappa:x^i_{\gamma_0} <_{\bool_i} x^i_{\gamma_1}\}\in D$ so again $y_{\gamma_0}<_\bool y_{\gamma_1}$. 

So far we have proved that $\Depth^+(\bool)\geq\lambda^+$. Likewise, $\Depth^+(\bool)\leq\lambda^+$ (when $\lambda^\kappa=\lambda$, hence $|\bool|=\lambda$) so part $(\beth)$ is established.
By claim \ref{cccc} below we shall get $\Depth^+(\bool_i)=\lambda$ for every $i<\kappa$, so the proof is accomplished.

\hfill \qedref{mt}

\begin{claim}
\label{cccc} Low $\Depth^+$ for every $\bool_i$. \newline 
$\Depth^+(\bool_i)=\lambda$ for every $i<\kappa$ in the construction above.
\end{claim}

\par\noindent\emph{Proof}.\newline 
Let $\langle\theta_\varepsilon: \varepsilon<\partial\rangle$ be as in the proof above, and let $\xi(\alpha)$ be the block function defined in that proof.
Fix any ordinal $i<\kappa$. For every $\varepsilon<\partial$, the sequence $\langle x^i_\alpha:\alpha\in[\theta_\varepsilon,\theta_{\varepsilon+1})\rangle$ is an increasing sequence in $\bool_i$, hence $\theta_{\varepsilon+1}< \Depth^+(\bool_i)$ for every $\varepsilon<\partial$. It means that $\lambda={\rm sup}\{\theta_{\varepsilon+1}:\varepsilon<\partial\}\leq\Depth^+(\bool_i)$.

Assume towards contradiction that $\bar{b}=\langle b_\gamma:\gamma<\lambda\rangle$ is an increasing sequence in $\bool_i$. Every member $b_\gamma\in\bool_i$ can be described by a Boolean term and a finite set of generators, $b_\gamma= \sigma_\gamma(\ldots,x^i_{\alpha(\gamma,\ell)},\ldots)_{\ell<n(\gamma)}$. Since $\cf(\lambda)=\partial>\aleph_0$ and there are just $\aleph_0$-many Boolean terms, we can assume without loss of generality that every $b_\gamma$ in our increasing sequence is generated by the same term $\sigma$ (in particular, there exists a natural number $n$ so that $n(\gamma)=n$ for every $\gamma<\lambda$). So we may write:

$$
b_\gamma=\sigma(\ldots,x^i_{\alpha(\gamma,\ell)},\ldots)_{\ell<n}
$$

We may assume (without loss of generality) that the finite sequence $\langle\alpha(\gamma,\ell):\ell<n\rangle$ is an increasing sequence of ordinals (for every $b_\gamma$). Observe that each ordinal $\alpha(\gamma,\ell)$ lies in a unique interval $[\theta_{\zeta(\gamma,\ell)}, \theta_{\zeta(\gamma,\ell)+1})$, which means that $\xi(\alpha(\gamma,\ell))=\zeta(\gamma,\ell)$.

By Lemma \ref{dddd} we can make (without loss of generality) the following assumptions. We assume that $T=\partial$ in the lemma, so for every $\varepsilon<\partial$ we have a finite set $r_\varepsilon\subseteq n$, acting as the root of the collection $\{\langle\alpha(\gamma,\ell): \ell<n\rangle: \gamma\in[\theta_\varepsilon, \theta_{\varepsilon+1})\}$. It means that the intersection of $\{\alpha(\gamma_0,\ell):\ell<n\}$ and $\{\alpha(\gamma_1,\ell):\ell<n\}$ equals $\{\alpha^\varepsilon_\ell: \ell\in r_\varepsilon\}$ for every distinct $\gamma_0,\gamma_1\in[\theta_\varepsilon, \theta_{\varepsilon+1})$.

Likewise, we assume that the collection $\{r_\varepsilon:\varepsilon<\partial\}$ is a Delta-system whose root is $r_*$. It means that $r_{\varepsilon_0}\cap r_{\varepsilon_1}=r_*$ for every $\varepsilon_0<\varepsilon_1<\partial$. Finally, if $\gamma_0\in[\theta_\varepsilon, \theta_{\varepsilon+1})$ and $\gamma_1\in[\theta_\zeta, \theta_{\zeta+1})$ then the intersection of $\{\alpha(\gamma_0,\ell):\ell<n\}$ and $\{\alpha(\gamma_1,\ell):\ell<n\}$ equals $\{\alpha_\ell:\ell\in r_*\}$.

The following property is important for the arguments below:

\begin{center}
$(\ast)$ We may assume that the finite sequence $\langle g_i(\zeta(\gamma,\ell)):\ell< n\rangle$ \newline 
does not depend on $\gamma$. \newline 
Namely, this is the same sequnece of ordinals \newline 
for every $\gamma<\lambda$.
\end{center}

Let us explain why this assumption can be made. The ordinal $\gamma_i$ from Lemma \ref{tttransitivity} is less than $\kappa^+$, and we have but $\kappa$-many $\gamma_i$-s. Hence $\delta = \sup\{\gamma_i:i<\kappa\}<\kappa^+ \leq\partial = \cf(\partial)$.

Each sequence of the form $\langle g_i(\zeta(\gamma,\ell)):\ell< n\rangle$ is an element of $[\delta]^{<\omega}$, so the number of possible sequences is strictly less than $\partial$. Since $\lambda>\cf(\lambda)=\partial$ we may assume that we have the same sequence for every $\gamma<\lambda$.

We may assume, in addition, that for some $S\subseteq\partial,|S|=\partial$ we have the following:

$$
\varepsilon_0,\varepsilon_1\in S, \varepsilon_0<\varepsilon_1 \Rightarrow \bigwedge\limits_{\gamma\in[\theta_{\varepsilon_0}, \theta_{\varepsilon_0+1})} \bigwedge\limits_{\ell<n} \zeta(\gamma,\ell)<\varepsilon_1.
$$

Actually, the set $S$ for which the proviso above is satisfied is a club subset of $\partial$. 
Fix two ordinals $\varepsilon_1,\varepsilon_2\in S$, such that $\varepsilon_1<\varepsilon_2$. Choose any $\gamma_1\in [\theta_{\varepsilon_1},\theta_{\varepsilon_1+1})$ and $\gamma_2\in [\theta_{\varepsilon_2},\theta_{\varepsilon_2+1})$. Set:

$$
Y=\{x^i_{\alpha(\gamma_1,\ell)}:\ell<n\} \bigcup \{x^i_{\alpha(\gamma_2,\ell)}:\ell<n\}.
$$

Set $\Gamma'_i=\{\varphi\in\Gamma_i:\varphi$ mentions only members of $Y\}$.
Let $\bool_Y$ be the Boolean algebra generated freely from the members of $Y$, except the relations mentioned in $\Gamma'_i$. Without loss of generality, $\bool_Y\subseteq\bool_i$ (for this, see \cite{MR1812172}, \S 3). Since $\gamma_1<\gamma_2<\lambda$, $\bool_i \models b_{\gamma_1}<b_{\gamma_2}$. As all the generators mentioned in $b_{\gamma_1},b_{\gamma_2}$ belong to $Y$ we have $\bool_Y \models b_{\gamma_1}<b_{\gamma_2}$ as well.

We define a function $f:Y\rightarrow Y$ as follows.
For every $\ell<n$ we define:

$$
f(x^i_{\alpha(\gamma_1,\ell)})=x^i_{\alpha(\gamma_2,\ell)} \qquad
f(x^i_{\alpha(\gamma_2,\ell)})=x^i_{\alpha(\gamma_1,\ell)}
$$

Notice that $f$ is a well-defined permutation of $Y$ (by the Delta-system requirements) of order $2$, i.e., $f\circ f={\rm Id}_Y$.
We claim that $f$ maps $\Gamma'_i$ onto itself. 
Let us prove this statement.

A typical member of $\Gamma'_i$ is an inequality $\eta=(x^i_{\alpha(\gamma_{j_1},\ell_1)}\leq x^i_{\alpha(\gamma_{j_2},\ell_2)})$ when $j_1,j_2\in\{1,2\}, \ell_1,\ell_2\in n$ and $f(\eta)$ is the inequality $(f(x^i_{\alpha(\gamma_{j_1},\ell_1)})\leq f(x^i_{\alpha(\gamma_{j_2},\ell_2)}))$. Our goal is to show that $\eta\in\Gamma'_i$ iff $f(\eta)\in\Gamma'_i$.
For proving this, we distinguish five cases: \newline 

\emph{Case 1}: $\ell_1,\ell_2\in r_*$.

We shall prove that $(x^i_{\alpha(\gamma_{j_1},\ell_1)}\leq_{\bool_i} x^i_{\alpha(\gamma_{j_2},\ell_2)})\in \Gamma'_i \Leftrightarrow (x^i_{\alpha(\gamma_{3-j_1},\ell_1)}\leq_{\bool_i} x^i_{\alpha(\gamma_{3-j_2},\ell_2)})\in\Gamma'_i$.
Under the assumption $\ell_1,\ell_2\in r_*$ we have $x^i_{\alpha(\gamma_1,\ell_1)}=x^i_{\alpha(\gamma_2,\ell_1)}$ and $x^i_{\alpha(\gamma_1,\ell_2)}=x^i_{\alpha(\gamma_2,\ell_2)}$. It means that the inequality after applying $f$ is just the same. \newline 

\emph{Case 2}: $\ell_1,\ell_2\notin r_*$, and $\ell_1=\ell_2$. 

Let $\ell$ denote the common value of $\ell_1,\ell_2$. If $j_1=j_2$ then the inequality $x^i_{\alpha(\gamma_{j_1},\ell)}\leq x^i_{\alpha(\gamma_{j_2},\ell)}$ is just an identity, and trivially preserved under $f$. If $j_1\neq j_2$ then (since $\ell\notin r_*$) we have $\neg[\xi(\alpha(\gamma_{j_1},\ell))<_i\xi(\alpha(\gamma_{j_2},\ell))]$ so the inequalities $x^i_{\alpha(\gamma_{j_1},\ell)}\leq x^i_{\alpha(\gamma_{j_2},\ell)}$ and $x^i_{\alpha(\gamma_{j_2},\ell)}\leq x^i_{\alpha(\gamma_{j_1},\ell)}$ do not belong to $\Gamma_i$ (and consequently, not to $\Gamma'_i$).

The above cases cover all the possibilities of $\ell_1=\ell_2$, so without loss of generality $\ell_1\neq\ell_2$ and at least one of them does not belong to $r_*$. \newline 

\emph{Case 3}: $j_1=j_2$.

We have to show that $(x^i_{\alpha(\gamma_1,\ell_1)}\leq_{\bool_i} x^i_{\alpha(\gamma_1,\ell_2)})\in \Gamma'_i$ iff $(x^i_{\alpha(\gamma_2,\ell_1)}\leq_{\bool_i} x^i_{\alpha(\gamma_2,\ell_2)})\in\Gamma'_i$.
This holds by the properties of the Delta-system and the property $(\ast)$ above. \newline 

\emph{Case 4}: $j_1\neq j_2$, and $\zeta(\gamma_{j_1},\ell_1)\neq \zeta(\gamma_{j_2},\ell_2)$.

By symmetry, without loss of generality $j_1=1$ and $j_2=2$. Also, we may assume that $\ell_1<\ell_2$. 
From $(\ast)$ we know that $g_i(\zeta(\gamma_1,\ell_1)) = g_i(\zeta(\gamma_2,\ell_1))$ and $g_i(\zeta(\gamma_2,\ell_2)) = g_i(\zeta(\gamma_1,\ell_2))$. Hence $g_i(\zeta(\gamma_1,\ell_1)) < g_i(\zeta(\gamma_2,\ell_2))$ iff $g_i(\zeta(\gamma_2,\ell_1)) < g_i(\zeta(\gamma_1,\ell_2))$. In the language of $<_i$ we can write $\zeta(\gamma_1,\ell_1) <_i \zeta(\gamma_2,\ell_2)$ iff $\zeta(\gamma_2,\ell_1) <_i \zeta(\gamma_1,\ell_2)$, see Lemma \ref{tttransitivity}(b). But this means that $x^i_{\alpha(\gamma_1,\ell_1)}<_{\bool_i}x^i_{\alpha(\gamma_2,\ell_2)}$ iff $x^i_{\alpha(\gamma_2,\ell_1)}<_{\bool_i}x^i_{\alpha(\gamma_1,\ell_2)}$, as required. \newline 

\emph{Case 5}: $j_1\neq j_2$, and $\zeta(\gamma_{j_1},\ell_1)= \zeta(\gamma_{j_2},\ell_2)$.

This case follows from the Delta-system properties and $(\ast)$. \newline 

With $f$ at hand, we employ Theorem \ref{ssikorski} which ensures the existence of a Boolean automorphism $\hat{f}:\bool_Y\rightarrow\bool_Y$ extending $f$. It follows that $\hat{f}(b_{\gamma_1})=b_{\gamma_2}$ and $\hat{f}(b_{\gamma_2})=b_{\gamma_1}$, contradicting the order preservation property of any Boolean homomorphism.

\hfill \qedref{cccc}

The above construction works equally well while replacing well-ordered sets by linearly-ordered sets. This yields the following:

\begin{corollary}
\label{llll} A Length gap. \newline 
Assume $\lambda>\cf(\lambda)=\partial$, and $\partial\in[\kappa^+,2^\kappa]$. \newline 
Then we can find $D$ and $\langle\bool_i:i<\kappa\rangle$ such that:
\begin{enumerate}
\item [$(\aleph)$] $\Length^+(\bool_i)\leq\lambda$ for every $i<\kappa$.
\item [$(\beth)$] $\Length^+(\bool)\geq\lambda^+$ (and equality holds if $\lambda^\kappa=\lambda$).
\end{enumerate}
Consequently, $\prod\limits_{i<\kappa}\Length^+(\bool_i)/D< \Length^+(\prod\limits_{i<\kappa}\bool_i/D)$.
\end{corollary}

\par\noindent\emph{Proof}.\newline 
The same proof as above, upon noticing that we have used just the cardinality of the increasing sequence and not the well ordering of it.

\hfill \qedref{llll}

\begin{remark}
\label{eeee} 
\begin{enumerate}
\item [$(\alpha)$] It seems that the assumption $\lambda^\kappa=\lambda$ (for both theorems, about $\Depth^+$ and $\Length^+$) can be weakened. Anyway, some assumption of this kind is needed, as if $2^\kappa>\lambda$ then the theorems may fail (unless we add further assumptions).
\item [$(\beta)$] By Theorem \ref{xxxx}, $D$ can be chosen as a regular ultrafilter. Nonetheless, it seems that the existence of a $\partial$-increasing chain in $\kappa^\kappa/D$ is essential (and we hope to prove it elsewhere).
\end{enumerate}
\end{remark}

\newpage 

\bibliographystyle{amsplain}
\bibliography{arlist}

\end{document}